\documentclass[11pt]{article}
\usepackage{amsmath,amssymb,amsthm,amscd,esint}
\usepackage{url,endnotes,hyperref}
\usepackage{tikz-cd}

\numberwithin{equation}{section}

\newtheorem{theorem}{Theorem}[section]
\newtheorem{lemma}[theorem]{Lemma}
\newtheorem{definition}[theorem]{Definition}

\theoremstyle{corollary}
\newtheorem{corollary}[theorem]{Corollary}
\theoremstyle{conjecture}

\newtheorem{question}[theorem]{Question}

\theoremstyle{assumption}

\theoremstyle{proposition}
\newtheorem{proposition}[theorem]{Proposition}
\theoremstyle{remark}
\newtheorem{remark}[theorem]{Remark}
\numberwithin{equation}{section}
\everymath{\displaystyle}

\newcommand{\Ric}{\operatorname{Ric}}
\newcommand{\Hess}{\nabla^2}

\newcommand{\tr}{\operatorname{tr}}

\newcommand{\eps}{\varepsilon}
\newcommand{\C}{\mathbb C}
\newcommand{\R}{\mathbb R}

\newcommand{\ddc}{\sqrt{-1}\partial \bar{\partial}}

\newcommand{\D}{\sqrt{-1}\,\partial\bar\partial}
\newcommand{\reg}{\mathrm{reg}}

\begin{document}

\title{On the Structure of Kähler-Ricci Shrinkers}
\author{}
\date{}
\maketitle

\vspace{-3em}
\begin{center}
\Large
\begin{tabular}{cc}
Tongxin Xu & \qquad Zhenlei Zhang
\end{tabular}
\end{center}
\vspace{0.5em}

\begin{abstract}
We prove that every complete K\"ahler--Ricci shrinker with nonnegative
Ricci curvature is holomorphically isometric to
$\mathbb{C}^k\times N^{n-k}$ for some $0\le k\le n$, where
$\mathbb{C}^k$ is the Gaussian shrinker and $N$ is a compact
K\"ahler--Ricci shrinker with nonnegative Ricci curvature.
As a corollary, every complete K\"ahler--Ricci shrinker with
quasi-positive Ricci curvature is compact. We also prove that every singular K\"ahler--Ricci shrinker on a Fano cone is a Ricci-flat K\"ahler cone.
\end{abstract}
\maketitle
\tableofcontents

\section{Introduction}
The Ricci flow was introduced by Hamilton as an evolution equation for Riemannian metrics \cite{hamilton1982three}. In his study of self-similar solutions, Hamilton introduced Ricci solitons\cite{MR954419,MR1375255}, which evolve under the Ricci flow only by diffeomorphisms and scaling. Recall that a Ricci shrinker $(X,g,f)$ is a complete Riemannian manifold $(X,g)$ together with a smooth function $f$ satisfying
\begin{align*}
    \operatorname{Ric}+\nabla^2 f=\lambda g, \qquad \lambda>0.
\end{align*}



A particularly important class in real even dimensions is formed by Kähler-Ricci shrinkers. A Ricci shrinker $(X,g,f)$ is called a Kähler-Ricci shrinker if it admits a Kähler structure $J$ such that $\nabla f$ is a holomorphic vector field. In this case, one can rewrite the shrinker equation as
\begin{align*}
    \operatorname{Ric(\omega)}+\sqrt{-1}\partial\bar{\partial} f= \lambda \omega, \qquad \lambda>0,
\end{align*}
where $\omega$ is the Kähler form induced by $g$ and $J$, and $\operatorname{Ric}(\omega)$ its Ricci form. The additional complex structure brings methods from complex and algebraic geometry into the study of Kähler-Ricci shrinker.

A basic rigidity problem asks to what extent positivity of the Ricci tensor constrains the global geometry of a complete Ricci shrinker. In particular, the following question was recently emphasized in \cite{chow2025ricci}.

\begin{question}\label{que}
Must every complete Ricci shrinker satisfying $\operatorname{Ric}>0$ be compact?
\end{question}

This question has been answered affirmatively under several stronger curvature assumptions. In the Riemannian setting, Perelman showed that every three-dimensional noncollapsed shrinker with positive bounded sectional curvature is compact \cite{Perelman2003RicciFlowSurgery}. In arbitrary dimensions, Munteanu-Wang \cite{munteanu2017positively} proved that every complete Ricci shrinker with nonnegative sectional curvature and positive Ricci curvature is compact. Li-Ni \cite{LiNi2020Orthogonal} subsequently established compactness under the weaker assumption of weakly PIC$_1$ curvature together with positive Ricci curvature. Further related results can be found in \cite{NiWallach2008Classification,PetersenWylie2010Classification,naber2007noncompact,LiNiWang2018PIC,wu2023shrinkers,QuWu2024Compact}, among others.

In the K\"ahler setting, Ni \cite{ni2005ancient} proved compactness under nonnegative holomorphic bisectional curvature and positive Ricci curvature; an alternative proof was later given by Wu-Zhang \cite{wu2016remarks}. Zhang \cite{zhang2019gradient} and Li-Ni \cite{LiNi2020Orthogonal}, extended this result to nonnegative orthogonal bisectional curvature. In complex dimension two, compactness follows
from the classification of K\"ahler--Ricci shrinkers
\cite{conlon2024classification,cifarelli2024finite,bamler2024new,li2026Kähler}.
For direct proofs in dimension two and in the toric setting, see
\cite{xu2026compactness}.

In this paper, we answer Question \ref{que} for any Kähler-Ricci shrinker without symmetric assumption and other curvature conditions. In fact, we characterize all Kähler–Ricci shrinkers with nonnegative Ricci curvature.

\begin{theorem}\label{main}
Let $(X^n,g,J,f)$ be a complete K\"ahler--Ricci shrinker with
$\Ric_g\ge0$. Then $X$ is holomorphically isometric to
$\mathbb C^k\times N^{n-k}$ for some $0\le k\le n$, where
$\mathbb C^k$ is the Gaussian shrinker and $N$ is a compact
K\"ahler--Ricci shrinker with $\Ric\ge0$.
\end{theorem}

It would be interesting to know whether the factor $N$ in the splitting has positive Ricci curvature. We do not have an answer right now. Our proof combines the Fano fibration structure established by
Sun--Zhang \cite{sun2024k} with Berndtsson's positivity theorem for direct images \cite{berndtsson2009curvature}.
The former provides a proper holomorphic map with compact fibres
onto an affine variety. On the regular part of the base, the latter
gives nonnegative curvature for the pushforward of the K\"ahler
volume form which yields
monotonicity of the normalized volume densities under the induced
soliton flow. We then study the rescaled soliton flow and construct a semipositive
limit form that descends to the regular base. Combining normalized
fibre-volume convergence with the volume monotonicity yields an
exact volume identity. This identity and a Bochner argument produce
a parallel splitting of the tangent bundle with a flat factor.
The de Rham theorem gives the global holomorphic isometric splitting.

As a corollary, the Question \ref{que} is correct under the quasi-positive Ricci curvature. 
\begin{corollary}
Let $(X^n,g,J,f)$ be a complete Kähler--Ricci shrinker with $\Ric_g \ge0$, if $\Ric_g>0$ at some point $p$ of $X$, then $X$ is compact.
\end{corollary}

In fact, any Kähler–Ricci shrinker with constant scalar curvature has nonnegative Ricci curvature. As another consequence of our results, Theorem \ref{main} implies the following result
\begin{corollary}\label{rig}
Let \((X^n,g,J,f)\) be a complete K\"ahler--Ricci shrinker
with constant scalar curvature. Then \(X\) is rigid.
\end{corollary}

A Ricci shrinker is called rigid if it is isometric to a finite quotient of $N\times\mathbb{R}^k$, where $N$ is Einstein and $\mathbb{R}^k$ is Gaussian. Cao conjectured that every complete Ricci shrinker with constant scalar curvature is rigid. Petersen-Wylie \cite{PW09} and Fernández-López-García-Río \cite{FLGR16} established fundamental rigidity criteria. The conjecture was proved in dimension four by Cheng-Zhou \cite{CZ23}, for further results in higher dimensions, we refer to \cite{LOQW24,LOQW25,WW26,LQW26}. Very recently, corollary \ref{rig} was also established by Horácio-Li-Ou-Wu-Zhou \cite{HLOWZ26}.

We also study singular K\"ahler--Ricci shrinkers on Fano cones.
Ricci-flat K\"ahler cone metrics naturally satisfy the shrinker
equation. In the complete smooth case, the converse follows
from the uniqueness results of Li--Zhang \cite{LiZhang26}
and Conlon--Deruelle \cite{ConlonDeruelle26}, and was also
proved in our earlier work \cite{xu2026compactness}.
The following theorem extends this conclusion to the singular
setting and confirms a conjecture of Sun--Zhang
\cite[Conjecture~6.3]{sun2024k}.

\begin{theorem}\label{cone}
Let $(X,\xi)$ be a polarized Fano cone.
Every singular K\"ahler--Ricci shrinker on $(X,\xi)$
is a Ricci-flat K\"ahler cone metric.
\end{theorem}

\paragraph*{Acknowledgments.} The first author thanks Professor Mihai Păun for his lectures delivered at the Chinese Academy of Sciences in September 2026, from which he learned Berndtsson’s work. This inspired the present research.  Zhenlei Zhang is supported by the National Natural Science Foundation of China (Grant No. 12531001).

\paragraph*{Use of AI.} The author used ChatGPT Astra for research assistance, particularly in the development and proof of Lemma~3.4. The authors have independently verified all mathematical arguments and takes full responsibility for their correctness.

 
 


\section{Preliminaries}

\begin{definition}\label{soliton}
    A Ricci shrinker is a triple $(X,g,f)$, where $X$ is a complete Riemannian manifold with a vector field $\nabla^gf$ satisfying the equation
\[
\Ric(g) + \Hess_g(f) = g,
\]
the vector field $\nabla^g f$ is called the soliton vector field.

If $g$ is complete and K\"ahler with K\"ahler structure $J$, then we say that $(X,g,J,f)$ is a K\"ahler Ricci shrinker, $\nabla^g f$ is complete and real holomorphic, and
\begin{align}\label{eqq}
   \operatorname{Ric}(\omega) + \sqrt{-1}\partial \bar{\partial} f = \omega, 
\end{align}
where $\operatorname{Ric}(\omega) $ is the Ricci form of Kähler form $\omega$.
\end{definition}

The above definition implies the standard normalization conditions and soliton identities associated to the shrinker.

\begin{lemma}\label{iden}
    Let $(X, g, f)$ be a n-dimensional Ricci shrinker, with soliton vector field $\nabla^g f$ for a smooth real-valued function $f:M\to \mathbb{R}$. Then we have the following soliton identities:
    \begin{align*}
        R_g+|\nabla^g f|^2=2f, \qquad R_g+\Delta f=n, \qquad R\geq 0.
    \end{align*}
\end{lemma}
\begin{proof}
    The first two inequalities hold by direct computation,  See \cite{chow2023ricci} for the complete proof. $R\geq 0$ follows from \cite{chen2009strong}.
\end{proof}

The K\"ahler--Ricci shrinker admits richer algebraic-geometric structures.
 Here we recall the polarized Fano fibration
framework of \cite{sun2024k}. A \textit{fibration} is a surjective projective morphism $\pi \colon X \to Y $
between normal varieties such that $\pi_*\mathcal{O}_X = \mathcal{O}_Y$. 
\begin{definition}\cite{collins2019sasaki}
    A \textit{polarized affine cone} $(Y,\xi,T)$ consists of a normal affine variety $Y = \operatorname{Spec} R$ with a compact torus $T$-action admitting a unique fixed point, together with a vector $\xi \in \operatorname{Lie}(T)$ lying in the Reeb cone. Here the latter condition means that for the weight decomposition
\[
R = \bigoplus_{\alpha \in \operatorname{Lie}(T)^*} R_\alpha,
\]
one has $\langle \alpha, \xi \rangle > 0$ whenever $R_\alpha \neq 0$ and $\alpha \neq 0$.
\end{definition}
\begin{definition}[Polarized Fano fibration]
    A polarized Fano fibration $(\pi\colon X \to Y, \xi)$ is a fibration $\pi\colon X \to Y$ with the following properties:
\begin{enumerate}
    \item $\pi\colon X \to Y$ is a Fano fibration. That is, $X$ and $Y$ are normal varieties, $X$ is klt, and that $-K_X$ is $\pi$-ample and $\mathbb{Q}$-Cartier.
    \item $X$ and $Y$ are equipped with a $\pi$-equivariant torus action $T$, and $\xi \in \mathfrak{t} = \operatorname{Lie}(T)$.
    \item $(Y, T, \xi)$ is a polarized affine cone.
\end{enumerate}
\end{definition}

Sun-Zhang proved the following structural theorem via deep algebraic geometric techniques.
\begin{theorem}\cite[Theorem 3.1]{sun2024k}
    A K\"ahler-Ricci shrinker $(X,g,J,f)$ naturally defines a polarized Fano fibration. In particular, $X$ is quasi-projective.
\end{theorem}

Let $Y_{\mathrm{reg}}$ denote the smooth locus of normal affine variety $Y$. Let $\Delta\subset Y_{\mathrm{reg}}$ be the set of critical values of $\pi$, and $S:=Y_{\mathrm{reg}}\setminus\Delta$.
For smooth complete Kähler-Ricci shrinker $X$, the flow of $\nabla f$ also has equivariance. 

\begin{lemma}
\label{lem:flow-descent}
There is a unique one-parameter group of
biholomorphisms \(\Psi_s:Y\to Y\) satisfying $\pi\circ\Phi_s=\Psi_s\circ\pi$ for any $s\in\mathbb R$. Moreover,
\(S\) is invariant under
\(\Psi_s\). 
\end{lemma}

\begin{proof}
Let \(\Phi_s\) denote the flow generated by $\nabla f$. By the definition of $\pi:X\to Y$, we know \(\pi\) is the Remmert reduction of \(X\). By \cite[Section~7]{HeinznerIannuzzi1997},
each \(\Phi_s\) descends to a biholomorphism
\(\Psi_s:Y\to Y\) satisfying $\pi\circ\Phi_s=\Psi_s\circ\pi$. Since \(\pi\) is surjective, this descent is unique. The surjectivity also implies $\Psi_{s+t}=\Psi_s\circ\Psi_t$,
 $\Psi_0=\mathrm{id}_Y$ and
 $\Psi_s^{-1}=\Psi_{-s}$. The completeness of $g$ implies the completeness of \(\nabla f\) \cite{zhang2009completeness}, so the flow \(\Phi_s\) is defined for all \(s\in\R\). Each \(\Psi_s\) preserves \(Y_{\mathrm{reg}}\). $d\pi_{\Phi_s(x)}\circ d\Phi_s|_x
 =d\Psi_s|_{\pi(x)}\circ d\pi_x$ implies the rank
of \(d\pi\) is preserved, combining with $\Phi_s(F_y)=F_{\Psi_s(y)}$, we know \(\Delta\) and \(S\) are invariant.
\end{proof}

If $\pi:X \to Y$ is the identity map, then $(X,\xi)$ is called a \textit{Fano cone}. We also consider the following singular Kähler-Ricci shrinkers on Fano cones. 
\begin{definition}\cite[Definition 4.1]{SunZhangII}
    A singular Kähler--Ricci shrinker consists of the data \((X, \xi, \omega, f)\) satisfying the following conditions:
\begin{enumerate}
    \item \((X, \xi)\) admits a polarized Fano fibration structure;
    \item \(\omega\) is a \(\xi\)-invariant positive \((1, 1)\)-current with locally bounded potentials, whose restriction to \(X^{\mathrm{reg}}\) is a smooth Kähler form;
    \item \(f\) is a continuous function on \(X\), smooth on the regular locus \(X_{\mathrm{reg}}\). Moreover, \(f\) is proper and bounded from below, and satisfies \(\xi = J\nabla f\) on \(X_{\mathrm{reg}}\);
    \item \(\int_X |\nabla f|^k e^{-f} \omega^n < \infty\) for every \(k \geq 0\);
    \item on \(X_{\mathrm{reg}}\), the shrinker equation holds: $\operatorname{Ric}(\omega) + \sqrt{-1}\partial\bar{\partial}f = \omega$. 
\end{enumerate}
\end{definition}

\section{Kähler-Ricci shrinkers and Ricci curvature}
\subsection{The monotonicity of logarithmic volume density}
Let $\pi:X^{r+d}\to Y^d$ be a proper holomorphic submersion with compact
connected fibres from smooth complex manifold $X$ to $Y$. Let $L$ be a holomorphic hermitian line bundle over $X$. For the direct image bundle $E=\pi_*(K_{X/Y}\otimes L)$ with $L^2$ norm $\|\cdot\|$, Berndtsson established the following positivity result:
\begin{theorem}\cite[Theorem 1.2]{berndtsson2009curvature}\label{bee}
    If the total space \(X\) is K\"ahler and \(L\) is (semi)positive over \(X\), then \((E, \|\cdot\|)\) is (semi)positive in the sense of Nakano.
\end{theorem}
We consider $E$ as a line bundle, in this case, Nakano positivity is equivalent to the positivity of its Chern curvature form.
\begin{lemma}\label{lem:directimage}
    Let \(\pi:X^{r+d}\to Y^d\) be a proper surjective holomorphic
submersion with connected fibres and $(X,\omega)$ be a Kähler manifold with $\Ric(\omega)\geq0 $
. Then $\Theta_\mu \ge 0$ for the smooth positive volume form $\mu=\pi_*dV_{\omega}$
\end{lemma}
\begin{proof}
Equip \(L=K_X^{-1}\) with the Hermitian
metric \(h_L\) induced by \(\omega\). Its curvature is $\Theta(L,h_L)=\Ric(\omega)$. Properness and connectedness of the fibres imply
$\pi_*\mathcal O_X=\mathcal O_Y$. Hence,
\begin{align}\label{bundle}
     E=:\pi_*(K_{X/Y}\otimes K^{-1}_X)
   =\pi_*\pi^*(K^{-1}_Y)\simeq{K^{-1}_Y}\otimes \pi_*\mathcal{O}_X\simeq K^{-1}_Y.
\end{align}
In local submersion coordinates \((y,z)=(y_1,...,y_d,z_1,...,z_r)\), we set $G$ be the metric matrix on $X$ and $G_V$ be the matrix of induce metric on each fiber. We set
\begin{align*}
    e_X=\partial_{y_1}\wedge...\wedge \partial_{y_d}\wedge\partial_{z_1}\wedge...\wedge \partial_{z_r}. \qquad e_Y=\partial_{y_1}\wedge...\wedge \partial_{y_d}
\end{align*} 
\(K_{X/Y}\otimes K_X^{-1}\simeq\pi^*K_Y^{-1}\) implies
the local section \(dz\otimes e_X\) corresponds to \(\pi^*e_Y\).
Thus the frame \(e_Y\) of \(E\) is represented
on each fibre by \((dz\otimes e_X)|_{X_y}\).
Its fibrewise \(L^2\) norm is
\begin{align*}
    \|e_Y(y)\|^2_E=\int_{X_y}|e_Y|^2dV_{X_y}=\int_{X_y}\frac{\operatorname{det}G}{\operatorname{det}G_V}{\operatorname{det}G_V}d\lambda_z=\int_{X_y}{\operatorname{det}G}d\lambda_z:=m(y),
\end{align*}
where \(d\lambda_z\) denotes the Euclidean volume form. Moreover, Fubini theorem implies $\mu=m(y)d\lambda_y$.  Then theorem \ref{bee} implies
\begin{align*}
     \Theta(E,h_E)
=-\ddc\log h_E(e_Y,e_Y)
=-\ddc\log m
=\Theta_\mu\ge 0.
\end{align*}
\end{proof}
Let \(Y\) be a complex manifold with a smooth positive volume form
\(\mu\).
Let \(W\) be a complete real holomorphic vector field with flow
\(\Psi_s\). Fix \(a\in\mathbb{R}\) and a positive smooth reference volume form \(\Omega\), we consider the following logarithmic volume density
\begin{align*}
    h_p(s)=\log\frac{e^{-as}\Psi_s^*\mu}{\Omega}(p),
\end{align*}
we have the following concavity and monotonicity lemma.
\begin{lemma}\label{lem:concavity}
If \(\Theta_\mu\geq0\) and \(\mathcal L_{JW}\mu=0\), then 
\begin{equation}\label{eq:concavity}
 h_p''(s)=-2(\Theta_\mu)_{\Psi_s(p)}(W,JW)\leq0.
\end{equation}
If \(h_p\) is bounded below on \([0,\infty)\), then \(h_p'(s)\geq0\)
for any \(s\geq0\). 
\end{lemma}

\begin{proof}
Direct computation implies
\[
 h_p'(s)=-a+(\operatorname{div}_\mu W)(\Psi_s(p)),\qquad
 h_p''(s)=W(\operatorname{div}_\mu W)(\Psi_s(p)).
\]
If $W\neq 0$, in local holomorphic coordinate, we write \(z_1=x+\sqrt{-1}y\) and arrange
\(W=\partial_x\), \(JW=\partial_y\).
If \(\mu=m\,d\lambda\) in this holomorphic coordinate, $\mathcal L_{JW}\mu=0$ gives \(\partial_y m=0\), so
\[
 \operatorname{div}_\mu W=\partial_x\log m,\qquad
 \Theta_\mu(W,JW)=-\frac12\partial_x^2\log m.
\]
which proves \eqref{eq:concavity}.
If \(W(p)=0\), the formula also holds. Finally, if \(h_p'(s_0)<0\), concavity gives
\(h_p(s)\leq h_p(s_0)+(s-s_0)h_p'(s_0)\) for \(s\geq s_0\),
contradicting a lower bound. Hence \(h_p'\geq0\).
\end{proof}
Given a K\"ahler manifold $(X,\omega)$ equipped with a holomorphic fibration, we establish the following technical lemma.

\begin{lemma}\label{lem:volume}
Let \(\pi:X^{r+d}\to Y^d\) be a proper surjective holomorphic
submersion with connected fibres, and let \(B\Subset Y\)
be a coordinate ball. Let
\[
 \beta=\D b\geq0,\qquad b\in C^0(B,\R),\qquad
 \beta\in L^\infty_{\mathrm{loc}}(B).
\]
Let \(\omega_j=\eps_j\omega+\D U_j\) be smooth K\"ahler forms,
where \(\eps_j>0\), \(\eps_j\to0\), and \(U_j\in C^\infty(X,\R)\).
For any compact subsets $K$ of \(B\),
\[
 \pi^*\beta\leq\omega_j\leq C\omega,\qquad
 U_j \longrightarrow \pi^*b \quad\text{uniformly on } \pi^{-1}(K),
\]
with \(C\) independent of \(j\). Then
\begin{equation}\label{eq:volume-convergence}
 \eps_j^{-r}\pi_*\frac{\omega_j^{r+d}}{(r+d)!}
 \longrightarrow V_F\frac{\beta^d}{d!}
 \quad\text{in }L^1_{\mathrm{loc}}(B),
 \qquad V_F:=\int_{X_y}\frac{\omega^r}{r!}.
\end{equation}
\end{lemma}
\begin{proof}
Set
\[
 u_j:=U_j-\pi^*b \qquad
 \eta_j=\omega_j-\pi^*\beta
       =\eps_j\omega +\D u_j\geq0.
\]
On each compact fibre \(X_y\), the restriction
\(\eta_j|_{X_y}=\omega_j|_{X_y}\) is smooth and represents
\(\eps_j[\omega|_{X_y}]\). Hence
\begin{equation}\label{eq:fibre-cohomology}
 \int_{X_y}\eta_j^{r-k}\wedge(\eps_j\omega)^k
 =\eps_j^r r!V_F,\qquad 0\leq k\leq r.
\end{equation}

The differential identities below are justified by local convolution
\(b_\delta=b*\rho_\delta\), for which
\[
 b_\delta\to b \quad\text{in }C^0_{\mathrm{loc}}(B),\qquad
 \D b_\delta\to\beta \quad\text{a.e.},\qquad
 \sup \nolimits_{\delta\ll1}\|\D b_\delta\|_{L^\infty(K)}<\infty
 \quad \text{on } K\Subset B.
\]
We omit $\delta$ from the notation and pass to the limit as currents in what follows.
Since \(\eta_j-\eps_j\omega=\D u_j\), direct computations implies
\[
 \eta_j^{r+1}-(\eps_j\omega)^{r+1}
 =\D\left(
 u_j\sum_{k=0}^r\eta_j^{r-k}\wedge(\eps_j\omega)^k
 \right).
\]
Then we have
\begin{align}
 T_j:=\eps_j^{-r}\pi_*(\eta_j^{r+1})
 =\eps_j\pi_*(\omega^{r+1})+\D p_j,
 \label{eq:transgression} \quad
 p_j=\eps_j^{-r}\pi_*\left(
 u_j\sum_{k=0}^r\eta_j^{r-k}\wedge(\eps_j\omega)^k
 \right).
 \notag
\end{align}
By \eqref{eq:fibre-cohomology},
\[
 \|p_j\|_{L^\infty(K')}
 \leq (r+1)r!V_F
       \|u_j\|_{L^\infty(\pi^{-1}(K'))}.
\]

Fix a flat K\"ahler form \(\omega_{B}\) on \(B\).
Given \(K\Subset\operatorname{int}K'\Subset B\), choose
\(\chi\in C_c^\infty(B)\) with
\(0\leq\chi\leq1\), \(\chi=1\) near \(K\), and
\(\operatorname{supp}\chi\subset K'\).
Since \(T_j\geq0\), integration by parts gives
\begin{equation}\label{eq:mass-estimate}
\begin{aligned}
 0\leq\int_K T_j\wedge\omega_{B}^{d-1}
 \leq\int_B\chi T_j\wedge\omega_{B}^{d-1}
 &=\eps_j\int_B\chi\,\pi_*(\omega^{r+1})
                       \wedge\omega_{B}^{d-1}
   +\int_B p_j\,\D\chi\wedge\omega_{B}^{d-1}\\
 &\leq C_K\left(
 \eps_j+\|u_j\|_{L^\infty(\pi^{-1}(K'))}
 \right)\longrightarrow0.
\end{aligned}
\end{equation}

To control the higher powers of $\eta_j$, we express it in matrix form with respect to a fixed smooth horizontal splitting over $\pi^{-1}(K')$:
\[
 \eta_j=
 \begin{pmatrix}
  H_j&B_j\\
  B_j^*&D_j
 \end{pmatrix},
 \qquad
 S_j=H_j-B_jD_j^{-1}B_j^*.
\]
Here \(D_j>0\) is the vertical fiber block. Since
\(0\leq\eta_j\leq C\omega\),
\[
 0\leq S_j\leq H_j\leq C_K\omega_{B}.
\]
The block determinant formula shows that the coefficient of $\eta_j^{r+k}\wedge\pi^*\omega_B^{d-k}$ relative to a fixed reference volume form is given by $\det D_j\,\sigma_k(S_j)$, where \(\sigma_k(S_j)\) denotes the \(k\)-th elementary symmetric
polynomial of the eigenvalues of \(S_j\) relative to \(\omega_B\).
The eigenvalues of \(S_j\) relative to \(\omega_{B}\) are uniformly
bounded and nonnegative, so
\[
 \sigma_k(S_j)\leq C_{K,k}\tr_{\omega_{B}} S_j.
\]
Consequently,
\begin{equation}\label{eq:higher-powers}
 \eta_j^{r+k}\wedge\pi^*\omega_{B}^{d-k}
 \leq C_{K,k}\,
       \eta_j^{r+1}\wedge\pi^*\omega_{B}^{d-1},
 \qquad 1\leq k\leq d.
\end{equation}
Finally, applying
\eqref{eq:fibre-cohomology}, we obtain
\begin{equation}\label{eq:volume-expansion}
 \eps_j^{-r}\pi_*\frac{\omega_j^{r+d}}{(r+d)!}
 =V_F\frac{\beta^d}{d!}+R_j,
 \qquad
 R_j:=
 \sum_{k=1}^d
 \frac{\eps_j^{-r}\pi_*(\eta_j^{r+k})\wedge\beta^{d-k}}
      {(r+k)!(d-k)!}\geq0.
\end{equation}
Using \(\beta\leq C_K\omega_{B}\),
\eqref{eq:mass-estimate}, and \eqref{eq:higher-powers}, we conclude
\[
 \|R_j\|_{L^1(K)}
 =\int_K R_j
 \leq C_K\left(
 \eps_j+\|u_j\|_{L^\infty(\pi^{-1}(K'))}
 \right)\longrightarrow0.
\]
Together with \(R_j\geq0\), this proves both assertions.
\end{proof}

\subsection{Proof of Theorem \ref{main}}
We first fix some setup. Let $X$ be a complete Kähler-Ricci shrinker with $\Ric \ge 0$ and $\pi:X^{r+d}\to Y^d$ is the associated polarized Fano fibration. Let 
\[
 S:=\left\{
 y\in Y_{\mathrm{reg}}:
 \operatorname{rank}_{\C}d\pi_x=d
 \text{ for every }x\in\pi^{-1}(y)
 \right\},
 \qquad X^\circ:=\pi^{-1}(S).
\]
Then \(\pi|_{X^\circ}:X^\circ\to S\) is a proper
holomorphic submersion.

Let $\Phi_s$ be the flow generated by $\nabla f$, and let $\Psi_s$ denote the flow obtained from Lemma \ref{lem:flow-descent}, where $W$ is the corresponding vector field.
Set $\omega_s=e^{-2s}\Phi_s^*\omega$ and consider the following measures on $S$
\[
 \mu=\pi_*\frac{\omega^n}{n!},\qquad
 \mu_s=e^{-2rs}\pi_*\frac{\omega_s^n}{n!}, \qquad \nu:=V_F\frac{\beta^d}{d!}
\]
We first implies following volume identity.

\begin{proposition}\label{eqqq}
    For any $s\ge 0$, $\mu_s=\nu=\mu$ as measures.
\end{proposition}
\begin{proof}
     For the flow \(\Phi_s\) of \(\nabla f\), \eqref{eqq} implies 
\begin{equation}\label{eq:metric-monotonicity}
 \partial_s\omega_s=-2e^{-2s}\Phi_s^* \Ric\leq0,
 \qquad 0<\omega_s\leq\omega.
\end{equation}
Define the following functions on \(X\):
\[
 F_s=e^{-2s}f\circ\Phi_s,\qquad
 U_s=2e^{-2s}\int_0^s f\circ\Phi_t\,dt.
\]
Direct differentiation gives
\begin{equation}\label{eq:potential-estimates}
 \partial_s F_s=-e^{-2s} R\circ\Phi_s\leq0,
 \qquad
 |dF_s|_g^2\leq |dF_s|_{g_s}^2
 =e^{-2s}|\nabla f|^2\circ\Phi_s\leq2F_s\leq2f.
\end{equation}
Thus \(F_s\) decreases locally uniformly to a  locally
Lipschitz function \(G\). Moreover,
\begin{equation}\label{eq:potentials}
 U_s=2\int_0^s e^{-2(s-t)}F_t\,dt\longrightarrow G,
 \qquad \omega_s=e^{-2s} \omega+\D U_s,
\end{equation}
locally uniformly for the functions and distributionally for their
complex Hessians. The metric identity follows by integrating
\(\partial_s(\Phi_s^*\omega)=2\D(f\circ\Phi_s)\).
Consequently
\begin{equation}\label{eq:global-limit}
 \alpha:=\D G=\lim_{s\to\infty}\omega_s,
 \qquad 0\leq\alpha\leq\omega,
\end{equation}
where \(\alpha\) has locally bounded coefficients. The monotone
pointwise limit represents this current almost everywhere. On each compact connected fibre over \(S\), the distributional limit
of \(\omega_s|_{X_y}\) satisfies $\D(G|_{X_y})\geq0$. Thus \(G|_{X_y}\) is a plurisubharmonic function and is
constant. It follows that \(G=\pi^*b\) on \(X^\circ\), for a
continuous function \(b\) on \(S\). Put \(\beta=\D b\).
Since \(\pi\) is a holomorphic submersion
and \(\pi^*\beta=\alpha\) has locally bounded coefficients,
the same holds for \(\beta\).
On any coordinate ball \(B\Subset S\), \eqref{eq:global-limit} and \cite[Proposition 4.4]{LiZhang26} give
\begin{equation}\label{eq:base-lower-bound}
 \beta\geq c_K\omega_{B}>0\quad\text{almost everywhere on }K\Subset B,
\end{equation}
where \(\omega_{B}\) is a fixed flat K\"ahler form on \(B\). By monotonicity of $\omega_s$, \(\pi^*\beta=\alpha\leq\omega_s\).
On any
\(B\Subset S\), properness of $\pi$ and Lemma~\ref{lem:volume} imply
\begin{equation}\label{eq:volume-limit-main}
 \mu_s\geq\nu,\qquad
 \mu_s\longrightarrow\nu\quad\text{in }L^1_{\mathrm{loc}}(B).
\end{equation}
The lower bound holds for each \(s\geq0\) by 
\eqref{eq:volume-expansion}.
By \eqref{eq:base-lower-bound}, \(\nu\) is bounded below locally by
a fixed positive smooth volume form. For any \(s\), smoothness of
\(\mu_s\) extends this lower bound from almost everywhere to everywhere.

Lemma~\ref{lem:concavity} gives \(\Theta_\mu\geq0\).
The Killing field \(J\nabla f\) preserves \(dV_g\), the equivariance implies
\(\mathcal L_{JW}\mu=0\). Fibre integration also gives 
\begin{align*}
\mu_s
=e^{2rs}e^{-2ns}
  \pi_*\left(\Phi_s^*\frac{\omega^n}{n!}\right)
=e^{-2(n-r)s}\Psi_s^*
  \left(\pi_*\frac{\omega^n}{n!}\right)
=e^{-2ds}\Psi_s^*\mu.
\end{align*}
For fixed \(y\in S\), \eqref{eq:volume-limit-main} implies 
\(h_y(s)=\log(\mu_s/\mu)(y)\) bounded below on \([0,\infty)\).
Apply Lemma~\ref{lem:concavity} with \(a=2d\) and \(\Omega=\mu\), we know \(\mu_s(y)\) is nondecreasing. Fixing \(s_0\geq0\) and passing to the limit in \(\mu_s\geq\mu_{s_0}\) for \(s\geq s_0\), we have
\(\nu\geq\mu_{s_0}\) as measures. Combining with \eqref{eq:volume-limit-main}, we have
\begin{equation}\label{eq:volume-equality}
\ \mu_s=\nu=\mu\quad\text{for any } s\geq0.
\end{equation}
These are equalities of measures.
\end{proof}
The above proposition implies that $\nu$ admits a smooth representative $\mu$, which yields the regularity of $G$.
\begin{lemma}
    The limit function $G\in C^\infty(X)$ and satisfies $\Delta_g G = 2d$.
\end{lemma}
\begin{proof}
    On \(X^\circ\), the form \(\alpha=\pi^*\beta\) has complex rank
\(d\) almost everywhere. Its nonzero eigenvalues relative to \(\omega\)
satisfy \(0<\lambda_1,\ldots,\lambda_d\leq1\), and
\begin{equation}\label{eq:eigenvalue-volume}
 \frac{\alpha^d}{d!}\wedge\frac{\omega^r}{r!}
 =\left(\prod_{i=1}^d\lambda_i\right)\frac{\omega^n}{n!}.
\end{equation}
The projection formula and proposition \ref{eqqq} give
\[
 \pi_*\left(\frac{\alpha^d}{d!}\wedge\frac{\omega^r}{r!}\right)
 =V_F\frac{\beta^d}{d!}=\pi_*\frac{\omega^n}{n!}.
\]
Thus the nonnegative measure
\(\bigl(1-\Pi_i\lambda_i\bigr)dV_g\) has zero pushforward and is
zero. Consequently every \(\lambda_i=1\) almost everywhere, so we have
\begin{equation}\label{eq:laplacian}
 \Delta_gG=2\tr_\omega\alpha=2d
\end{equation}
almost everywhere on \(X^\circ\). \(\alpha=\D G\) has locally bounded coefficients, so its
trace has no singular part on \(X\setminus X^\circ\).
Since \(X\setminus X^\circ\) has measure zero, \eqref{eq:laplacian} holds
distributionally on \(X\). Elliptic regularity gives
\(G\in C^\infty(X)\).
\end{proof}

In fact, the smoothness of $G$ implies the smoothness of $\beta$, so we can obtain the following Riemannian submersion.
\begin{lemma}
    $\pi:(X^\circ,g)\longrightarrow(S,g_\beta)$ is a Riemannian submersion.
\end{lemma}
\begin{proof}
    Define the smooth endomorphism $P$ by
\begin{equation}\label{eq:projection}
 \alpha(U,V)=g(JPU,V),\qquad U,V\in TX.
\end{equation}
Since $\alpha$ is a real $(1,1)$-form, $P$ is self-adjoint
and commutes with $J$. 
We note that the eigenvalues of $\alpha$ correspond to the eigenvalues of $P$, so
\begin{align}\label{proj}
    P^2=P,\qquad \tr_{\R}P=2d\quad\text{on }X, 
\end{align}
We set \(H=\operatorname{Im}P\) and \(E=\ker P\). Since \(G=\pi^*b\) is smooth and \(\pi:X^\circ\to S\) is a
holomorphic submersion, \(b\) is smooth.
Thus \(\beta=\ddc b\) is smooth, we know $\beta$ is a K\"ahler form on $S$ by
\eqref{eq:base-lower-bound}. The identity \(\alpha=\pi^*\beta\) now holds as an equality
of smooth forms. Since \(\beta\) is nondegenerate,
\[
 E=\ker P=\ker\alpha=\ker d\pi
 \qquad\text{on }X^\circ.
\]
Moreover, for horizontal vectors \(U,V\in H\),
\[
 g_\beta(d\pi U,d\pi V)
 =\alpha(U,JV)
 =g(PU,V)
 =g(U,V).
\]
where $g_{\beta}$ is the Riemannian metric induced by $\beta$. Hence $\pi:(X^\circ,g)\longrightarrow(S,g_\beta)$ is a Riemannian submersion with horizontal distribution \(H\).
\end{proof}
To characterize this splitting, we further compute several equations associated with $G$.
\begin{lemma}\label{ggg}
    The smooth limit function \(G\) satisfies $\Hess_gG(\cdot,\cdot)=g(P\cdot,\cdot)$. 
\end{lemma}
\begin{proof}
   Direct computations implies $dF_s=e^{-2s}\Phi_s^*(df)=g_s(\nabla f,\cdot)$. Using \(F_s\to G\) locally uniformly,
\(g_s\to g(P\cdot,\cdot)\) almost everywhere, and
\(0\le g_s\le g\), we pass to the limit in distributions
to obtain $\nabla G=P\nabla f$.

For fixed \(t\in\R\), letting \(s\to\infty\) in
\(F_s\circ\Phi_t=e^{2t}F_{s+t}\) gives
\(G\circ\Phi_t=e^{2t}G\), hence \((\nabla f)G=2G\).
Since \(P\) is an orthogonal projection,
\begin{equation}\label{eq:gradient-length}
 |\nabla G|^2=g(P\nabla f,P\nabla f)=g(P\nabla f,\nabla f)=(\nabla f)G=2G.
\end{equation}
\(\alpha=\D G\) and \eqref{eq:projection} imply
\[
 \frac12\bigl(\Hess_gG(U,V)+\Hess_gG(JU,JV)\bigr)
 =g(PU,V).
\]
Thus the \(J\)-invariant part of \(\Hess_gG\) is
\(g(P\cdot,\cdot)\). Orthogonality gives
\[
 |\Hess_gG|^2 =|g(P\cdot,\cdot)|^2+\bigl|\Hess_gG-g(P\cdot,\cdot)\bigr|^2
 =2d+\bigl|\Hess_gG-g(P\cdot,\cdot)\bigr|^2.
\]
Bochner's formula, together with \eqref{eq:laplacian}
and \eqref{eq:gradient-length}, now gives
\begin{align*}
 2d
 =\frac12\Delta_g|\nabla G|^2
 =|\Hess_gG|^2+\Ric_g(\nabla G,\nabla G)
 =2d+\bigl|\Hess_gG-g(P\cdot,\cdot)\bigr|^2
      +\Ric_g(\nabla G,\nabla G).
\end{align*}
\(\Ric_g\geq0\) implies $\Hess_gG=g(P\cdot,\cdot)$.
\end{proof}

Now we proof the Theorem \ref{main}.
\begin{proof}[Proof of Theorem \ref{main}]
    For the closed form \(\omega-\alpha\), the definition of $P$ gives
    \begin{align*}
\ker(\omega-\alpha)=\ker(I-P)=\operatorname{im}P=H.
    \end{align*}
For any vector field $W$ and $U,V\in H$, we have
\begin{align*}
0=d(\omega-\alpha)(U,V,W)=-(\omega-\alpha)([U,V],W),
\end{align*}
so $[U,V]\in H$ i.e. $H$ is integrable. On \(X^\circ\), the Riemannian submersion formula gives
\[
 (\nabla_U V)^{\mathrm{vert}}
 =\tfrac12[U,V]^{\mathrm{vert}}=0
\]
for horizontal fields \(U,V\in H \). Thus the horizontal leaves are
totally geodesic on \(X^\circ\).
Their second fundamental form is smooth and vanishes on the
dense set \(X^\circ\), so it vanishes on $X$, the horizontal leaves are totally geodesic on $X$. Combining with the completeness of $(X,g)$, we know each maximal horizontal leaf \(L\) is geodesically complete. By Lemma \ref{ggg}, we have $\Hess_L(G|_L)=g_L$, \cite[Theorem 5.75]{petersen2006riemannian} implies $(L,g_L)\cong\mathbb R^{2d}$.

Since the horizontal leaves are Euclidean and $\pi$ is a Riemannian submersion, we know 
the metric \(g_\beta\) is flat.
Choose holomorphic coordinates \(z_1,\ldots,z_d\)
on a neighborhood \(U_{S}\subset S\).
The functions $u_i=\operatorname{Re}(z_i\circ\pi)$ and $u_{d+i}=\operatorname{Im}(z_i\circ\pi)$ are harmonic on $\pi^{-1}(U_{S})$, and $|\nabla u_i|$ is constant since $\pi$ is a Riemannian submersion. Bochner's formula gives
\[
 0=\frac12\Delta_g|\nabla u_i|^2
  =|\Hess u_i|^2+\Ric_g(\nabla u_i,\nabla u_i).
\]
Both terms are nonnegative, so \(\nabla u_i\) is parallel.
These gradients span \(H\), hence 
\begin{align*}
P(\cdot)=\sum\nolimits_{i=1}^{2d}g(P(\cdot),\nabla u_i)\nabla u_i.
\end{align*}
So \(\nabla P=0\) on \(X^\circ\).
By smoothness and density, \(\nabla P=0\) on all of \(X\). Thus \(H\) and \(E\) are parallel distributions.
Since \(X\) is complete and simply connected
\cite[Proposition~3.10]{sun2024k}, the de Rham decomposition theorem gives $(X,g)\simeq(\R^{2d},g_{\mathrm E})\times(N,g_N)$, both distributions are \(J\)-invariant,
so this is a holomorphic isometric splitting
\[
 (X,g,J)\simeq
 (\C^d,g_{\mathrm E},J_{\mathrm E})\times(N,g_N,J_N).
\]

Choose a regular fibre \(F=\pi^{-1}(y)\), \(y\in S\).
It is compact and connected, with \(TF=E|_F\).
Its Euclidean coordinate is therefore constant, so
\(F\subset\{z_0\}\times N\) for some \(z_0\in\C^d\).
Since \(F\) and \(N\) have the same dimension and $F$ is compact, we know \(F=\{z_0\}\times N\). In particular, \(N\) is compact. The product metric and the shrinker equation imply that \(N\)
is a compact K\"ahler--Ricci shrinker with nonnegative Ricci
curvature. 
\end{proof}

\begin{proof}[Proof of Corollary \ref{rig}]
  We write \(R\equiv c\).
Let \(\Phi_s\) be the flow of \(\nabla f\).
Along the flow,
\[
 \frac{d}{ds}(f\circ\Phi_s)=2f\circ\Phi_s-c,
\qquad
 f\circ\Phi_s
 =\frac{c}{2}+e^{2s}\left(f-\frac{c}{2}\right).
\]
Then
\[
 \frac{d}{ds}\Phi_s^*\omega
 =2\D(f\circ\Phi_s)
 =2e^{2s}\D f.
\]
Integrating gives
\[
 \Phi_s^*\omega
 =\omega+(e^{2s}-1)\D f
 =\Ric + e^{2s}\D f.
\]
Letting \(s\to-\infty\), we obtain \(\Ric_g\geq0\). By Theorem \ref{main}, $X$ is biholomorphically isometric to $\mathbb{C}^{k}\times N^{n-k}$, where $N$ is a compact Kähler-Ricci shrinker with nonnegative Ricci curvature. Since $R\equiv c$ implies $R_N\equiv c$, we have $\Delta_N f_N = n-k-c$. Integrating implies $f_N$ is constant i.e. $N$ is a Kähler–Einstein manifold. So $X$ is rigid.
\end{proof}

\section{Singular Kähler-Ricci shrinkers on Fano cones}

In this section, let \((X,\xi,\omega,f)\) be a singular
K\"ahler--Ricci shrinker on a polarized Fano cone \((X,\xi)\).
Thus \(X\) is a normal affine variety with klt singularities,
equipped with an algebraic torus action and a positive Reeb
field \(\xi\). Denote its vertex by \(o\), its regular locus
by \(X_{\mathrm{reg}}\), and its singular set by
\(X_{\mathrm{sin}}\).

Write \(n=\dim_{\mathbb C}X\), and let \(g\) be the metric
induced by \(\omega\) on \(X_{\mathrm{reg}}\). Set
\[
 d\mu=e^{-f}\frac{\omega^n}{n!},\qquad
 \Delta_f=\Delta_g-\langle\nabla f,\nabla\cdot\rangle, \qquad \operatorname{div}_{\mu}(\cdot)= \operatorname{div}(\cdot) - \langle\nabla f,\cdot\rangle.
\]
By assumption, \(f\) is locally bounded, proper and bounded
below on \(X\), and \(0<\mu(X)<\infty\).
Choose \(m>0\) such that \(L=K_X^{-m}\) is a line bundle.
On \(X_{\mathrm{reg}}\), let \(h_g\) be the metric on \(L\)
induced by \(g\), and set \(H=e^{-mf}h_g\). The shrinker
equation gives
\[
 \frac{1}{m}\Theta_H
 =\operatorname{Ric}(\omega)
       +\sqrt{-1}\partial\bar\partial f
 =\omega.
\]
The normalized local weights of \(H\) differ from locally
bounded potentials of \(\omega\) by pluriharmonic functions.
By normality, these functions extend across the singular
set, hence \(H\) extends
to \(L\) with locally bounded weights. The algebraic torus action defines biholomorphisms
\[
 \Phi_t(x)=\exp(-it\xi)\cdot x,\qquad t\in\mathbb R.
\]
These preserve \(X_{\mathrm{reg}}\), where their generator
is \(\nabla f\). Thus \(\nabla f\) is complete on
\(X_{\mathrm{reg}}\), and its flow extends to
biholomorphisms of \(X\).

We first recall the construction of cut-off functions on singular spaces. The following lemma follows from \cite[Lemma 3.7]{SJ}.
\begin{lemma}\label{lem:singular-cutoffs}
For every compact set \(K\subset X\), there exist a neighborhood
\(U\) of \(K\) and functions \(\chi_j\in C^\infty(U\cap X_{\reg})\),
each vanishing near \(U\cap X_{\mathrm{sin}}\), such that
\[
 0\leq\chi_j\leq1,\qquad
 \chi_j\longrightarrow1
 \quad\text{locally uniformly on }U\cap X_{\reg},\qquad
 \int_{K\cap X_{\reg}}|\nabla\chi_j|_g^2\,d\mu\longrightarrow0.
\]
\end{lemma}

\begin{proof}
We take a resolution \(q:\widetilde{X}\to X\) and take \(K\subset U\Subset X\).
Let \(\sigma\) define an effective divisor supported on
\(E=q^{-1}(X_{\mathrm{sin}})\) and fix a smooth Hermitian metric $h$ such that
\begin{align*}
    \ell:=\log|\sigma|_{h}^2\leq0 \qquad \text{near} \quad q^{-1}(\overline U).
\end{align*}
Let \(F\) be a smooth convex nondecreasing function with $F=0$ on $(-\infty,0]$ and $F(1)=1$. Set $\chi_j = F(1+\ell/j)$. Each $\chi_j$ vanishes near $E$, extends smoothly by zero across $E$, and the sequence converges locally uniformly to $1$ away from $E$.
 Via the biholomorphism
\(q:q^{-1}(U)\setminus E\to U\cap X_{\mathrm{reg}}\),
we regard \(\chi_j\) as functions on \(U\cap X_{\mathrm{reg}}\).

Set \(T=(q^*\omega)^{n-1}\), and choose
\(\eta\in C_c^\infty(q^{-1}(U))\), with \(\eta\geq0\) and
\(\eta=1\) near \(q^{-1}(K)\).
Away from \(E\), convexity of \(F\) gives
\[
 \sqrt{-1}\partial\bar\partial\chi_j
 =\frac{F'(1+\ell/j)}{j}\,\sqrt{-1}\partial\bar\partial\ell
  +\frac{F''(1+\ell/j)}{j^2}\,
    \sqrt{-1}\partial\ell\wedge\bar\partial\ell
 \geq-\frac{C}{j}\theta,
\]
near \(\operatorname{supp}\eta\), where \(\theta\) is a fixed
smooth positive \((1,1)\)-form. Integration by parts yields
\[
 \begin{aligned}
 \int_{q^{-1}(U)}\eta\,\sqrt{-1}\partial\chi_j\wedge\bar\partial\chi_j\wedge T
 &=\frac12\int_{q^{-1}(U)}\chi_j^2 \sqrt{-1}\partial\bar\partial\eta\wedge T
   -\int_{q^{-1}(U)}\eta\chi_j \sqrt{-1}\partial\bar\partial\chi_j\wedge T\\
 &\leq\frac12\int_{q^{-1}(U)}\chi_j^2 \sqrt{-1}\partial\bar\partial\eta\wedge T
   +O(j^{-1}) \longrightarrow 0.
 \end{aligned}
\]
Combining with the local boundedness of \(f\) we have
\[
 \int_{K\cap X_{\reg}}|\nabla\chi_j|_g^2\,d\mu
 \longrightarrow0.
\]
\end{proof}

\begin{lemma}\label{lem:singular-integration}
Let \(v\in C^\infty(X_{\mathrm{reg}})\),
with \(v\) and \(\Delta_fv\) locally bounded on \(X\),
then

\begin{equation}\label{eq:singular-integration}
 \int_X\psi\,\Delta_fv\,d\mu
 =-\int_X\langle\nabla\psi,\nabla v\rangle\,d\mu
\end{equation}
for every bounded, compactly supported test function \(\psi\)
on \(X\) that is locally Lipschitz on \(X_{\mathrm{reg}}\)
and satisfies \(\int_X|\nabla\psi|^2\,d\mu<\infty\).
\end{lemma}

\begin{proof}

Fix a compact set $K\subset X$, and let
$q:\widetilde X\to X$ be a resolution.
Choose $\widetilde\eta\in C_c^\infty(\widetilde X)$ with
$0\leq\widetilde\eta\leq1$ and $\widetilde\eta=1$ near
$q^{-1}(K)$.
Define $\eta$ on $X_{\mathrm{reg}}$ by
$\eta\circ q=\widetilde\eta$ on $q^{-1}(X_{\mathrm{reg}})$.
Then $\eta$ is smooth, equals $1$ on $K\cap X_{\mathrm{reg}}$,
and its support has compact closure contained in
$q(\operatorname{supp}\widetilde\eta)$.
Since $f$ is locally bounded and $(q^*\omega)^{n-1}$
has locally finite mass,
\[
 \int_{X_{\mathrm{reg}}}|\nabla\eta|_g^2\,d\mu
 \leq C\int_{\widetilde X}
 \sqrt{-1}\partial\widetilde\eta
 \wedge\bar\partial\widetilde\eta
 \wedge(q^*\omega)^{n-1}<\infty.
\]
Apply Lemma~\ref{lem:singular-cutoffs} to the compact set
$q(\operatorname{supp}\widetilde\eta)$ to obtain $\chi_j$. Since \(\eta\chi_j\) has compact support in \(X_{\mathrm{reg}}\),
testing with \((\eta\chi_j)^2v\) gives
\begin{align*}
\int_X(\eta\chi_j)^2|\nabla v|^2\,d\mu
 &=-\int_X(\eta\chi_j)^2v\Delta_fv\,d\mu
 -2\int_X(\eta\chi_j)v
       \langle\nabla(\eta\chi_j),\nabla v\rangle\,d\mu\\
&
\leq
\int_X(\eta\chi_j)^2|v\Delta_fv|\,d\mu +\int_X (\frac12(\eta\chi_j)^2|\nabla v|^2
+2v^2|\nabla(\eta\chi_j)|^2)\,d\mu,
\end{align*}
since \(v\) and \(\Delta_fv\) locally bounded on \(X\), we have 
\[
 \int_X(\eta\chi_j)^2|\nabla v|^2\,d\mu
 \leq C\int_X
 \left((\eta\chi_j)^2+|\nabla(\eta\chi_j)|^2\right)\,d\mu,
\]
where \(C\) is independent of \(j\).
The right-hand side is uniformly bounded by the 
finiteness of \(\mu\), the finite energy of \(\eta\), and
Lemma~\ref{lem:singular-cutoffs}.
Fatou's lemma implies
\[
 \int_X|\nabla v|^2\,d\mu
 \leq\liminf_{j\to\infty}
 \int_X(\eta\chi_j)^2|\nabla v|^2\,d\mu<\infty.
\]

Choose
a compact set \(K\subset X\) containing the support of $\psi$.
Testing with \(\chi_j\psi\), using smooth approximation on
\(X_{\mathrm{reg}}\), gives
\[
 \int_X\chi_j\psi\,\Delta_fv\,d\mu
 =-\int_X\chi_j\langle\nabla\psi,\nabla v\rangle\,d\mu
  -\int_X\psi\langle\nabla\chi_j,\nabla v\rangle\,d\mu,
\]
we note
\[
 \left|
 \int_X\psi\langle\nabla\chi_j,\nabla v\rangle\,d\mu
 \right|
 \leq\|\psi\|_\infty
 \left(\int_K|\nabla\chi_j|^2\,d\mu\right)^{1/2}
 \left(\int_K|\nabla v|^2\,d\mu\right)^{1/2}
 \longrightarrow0.
\]
The other terms converge by dominated convergence, proving
\eqref{eq:singular-integration}.
\end{proof}

We next construct a global K\"ahler potential $\varphi$
and derive its relation to the soliton potential $f$.
\begin{lemma}
After adding a constant to \(f\), there exists a nonnegative
global K\"ahler potential \(\varphi\), smooth on \(X_{\reg}\),
such that
\[
 \sqrt{-1}\partial\bar\partial\varphi=\omega,\qquad
 \nabla f(\varphi)=2f,\qquad
 \Delta_f f=\kappa-2f,\qquad
 \Delta_f\varphi=2n-2f
\]
on \(X_{\reg}\), for some \(\kappa\in\mathbb R\).
\end{lemma}

\begin{proof}
Since $X$ is affine, $K_X^{-m}$ has a global algebraic
section that does not vanish at $o$.
Choose a torus weight component $s$ with $s(o)\neq0$.
Its zero set is closed and invariant under $\Phi_t$.
Since $\Phi_t(x)\to o$ as $t\to-\infty$ for every $x\in X$,
any zero of $s$ would force $s(o)=0$.
Thus $s$ is nowhere vanishing and defines a global frame. Let \(\Omega\) be the canonical measure determined by
\(s^{-1}\), normalized so that
\[
 |s|_{h_g}^{2/m}\Omega=dV_g,
\qquad
 \varphi=:-\frac1m\log|s|_H^2.
\]
Since \(H=e^{-mf}h_g\) and \(\Theta_H=m\omega\), we have
\begin{equation}\label{eq:global-potential}
  \sqrt{-1}\partial\bar\partial\varphi=\omega,\qquad
 d\mu=e^{-f}dV_g=|s|_H^{2/m}\Omega=e^{-\varphi}\Omega.
\end{equation}
Then \(\varphi\) is locally bounded and plurisubharmonic
on \(X\), and smooth on \(X_{\reg}\).
The weight property of \(s\) gives
\[
 \Phi_t^*\Omega=e^{\kappa t}\Omega
\]
for some \(\kappa\in\mathbb R\).
Both \(\mu\) and \(\Omega\) are invariant under \(\xi\),
so \(\xi\varphi=0\).
Cartan's formula yields
\[
 0=\mathcal L_\xi d^c\varphi=-2df+d(\nabla f(\varphi)),
\]
where \(d^c\varphi=-d\varphi\circ J\). Thus \(\nabla f(\varphi)-2f\) is constant on 
\(X_{\reg}\). After adding a constant to \(f\) and rescaling
\(H\), we have
\[
 \nabla f(\varphi)=2f.
\]
Since \(d\mu=e^{-\varphi}\Omega\) and $\mathcal L_{\nabla f}\mu=\operatorname{div}_\mu(\nabla f) \mu$, we obtain
\[
 \Delta_f f=\operatorname{div}_\mu (\nabla f)
 =\kappa-\nabla f(\varphi)=\kappa-2f.
\]
Moreover, \(\sqrt{-1}\partial\bar\partial\varphi=\omega\)
gives \(\Delta_g\varphi=2n\), hence
\[
 \Delta_f\varphi=2n-\nabla f(\varphi)=2n-2f.
\]

Along every nonconstant orbit, \(f(\Phi_t x)\) is strictly
increasing and has a finite backward limit by local
boundedness near \(o\). The identity \(\nabla f(\varphi)=2f\)
and boundedness of \(\varphi\) near \(o\) force this limit
to be zero. Thus \(f>0\) on \(X_{\reg}\setminus\{o\}\),
and \(f(o)=0\) if \(o\in X_{\reg}\).
Consequently, \(\varphi\) is nondecreasing along forward
orbits. Its local lower bound near \(o\) is 
a global lower bound. Multiplying \(s\) by a constant,
we may assume \(\varphi\geq0\).
\end{proof}
\begin{remark}\label{kkk}
    On $\mathbb C^n$, write $z_j=x_j+\sqrt{-1}y_j$. The standard
Reeb vector field $\xi$ is $\textstyle\sum_{j=1}^n
 (-y_j\partial_{x_j}+x_j\partial_{y_j})$, and the flow of $-J\xi$ is $\Phi_t(z)=e^tz$.
For the standard volume form
$\Omega=\wedge_{j=1}^n(dx_j\wedge dy_j)$, we have
\[
 \Phi_t^*\Omega
 =\wedge_{j=1}^n(e^t dx_j\wedge e^t dy_j)
 =e^{2nt}\Omega.
\]
Thus $\kappa=2n$ for the standard Reeb field.
\end{remark}
From the observation in Remark \ref{kkk}, we deduce the following lemma.
\begin{lemma}
With the normalization above, 
$\kappa=2n$ and \(\varphi-f\) is constant on \(X_{\reg}\). 
\end{lemma}

\begin{proof}
Choose a smooth nonincreasing function
\(\theta:[0,\infty)\to[0,1]\), equal to \(1\) on \([0,1]\)
and to \(0\) on \([2,\infty)\), and set
\(\theta_R=\theta(f/R)\).
Its support has compact closure in \(X\) by properness of \(f\).
The local energy bounds and integrations by parts below follow
from Lemma~\ref{lem:singular-integration}.

Testing \(\Delta_f\varphi=2n-2f\) with \(\theta_R\) gives
\[
 \int_X\theta_R f\,d\mu
 =n\int_X\theta_R\,d\mu
  +\frac1R\int_X f\theta'(f/R)\,d\mu
 \leq n\mu(X).
\]
The last term has absolute value at most
\(2\|\theta'\|_\infty\mu(\{f\geq R\})\), which tends to zero.
Thus monotone convergence gives
\begin{align*}
  \int_X f\,d\mu=n\mu(X).
\end{align*}
Testing the same equation with \(f\theta_R\) gives
\[
 \int_X\theta_R f^2\,d\mu
 =(n+1)\int_X\theta_R f\,d\mu
  +\frac1R\int_X f^2\theta'(f/R)\,d\mu
 \leq(n+1)\int_X f\,d\mu.
\]
Monotone convergence yields
\begin{equation}\label{eq:moments}
 \int_X f^2\,d\mu\leq n(n+1)\mu(X).
\end{equation}
Since \(\operatorname{div}_\mu (\nabla f)=\kappa-2f\), for any $t>0$ we have
\[
 \begin{aligned}
 0
 =\int_X(\Phi_t^*d\mu-d\mu)
 =\int_0^t\int_X\Phi_s^*((\kappa-2f)d\mu)\,ds
 =t\left(\kappa\mu(X)-2\int_X f\,d\mu\right).
 \end{aligned}
\]
Hence \(\kappa=2n\), and
\begin{equation}\label{eq:harmonic}
 \Delta_f f=\Delta_f\varphi=2n-2f,\qquad
 \Delta_f(\varphi-f)=0.
\end{equation}

For \(R\geq n\), test the first equation in
\eqref{eq:harmonic} with \((R-f)_+\).
Using \eqref{eq:moments}, we obtain
\[
 \begin{aligned}
 \int_{\{f<R\}}|\nabla f|^2\,d\mu
 &=2R\int_{\{f\geq R\}}(f-n)\,d\mu
   +2\int_{\{f<R\}}f(f-n)\,d\mu\\
 &\leq2\int_X f(f-n)\,d\mu
 \leq2n\mu(X).
 \end{aligned}
\]
Letting \(R\to\infty\) gives
\begin{equation}\label{eq:energy-f}
 \int_X|\nabla f|^2\,d\mu\leq2n\mu(X).
\end{equation}
For \(a>0\), test \(\Delta_f(\varphi-f)=0\) with
\(\theta_Re^{-a\varphi}\). Since
\(\nabla f(\varphi)=2f\), we have
\[
 \begin{aligned}
 a\int_X\theta_Re^{-a\varphi}|\nabla\varphi|^2\,d\mu
 ={}&2a\int_X\theta_R f e^{-a\varphi}\,d\mu
 +\frac1R\int_X\theta'(f/R)e^{-a\varphi}
       (2f-|\nabla f|^2)\,d\mu.
 \end{aligned}
\]
 Since
\(f,|\nabla f|^2\in L^1(\mu)\) and \(\varphi\geq0\), the last term tends to zero.
Letting \(R\to\infty\) and \(a\to0\),
monotone convergence gives
\[
 \int_X|\nabla\varphi|^2\,d\mu
 =2\int_X f\,d\mu=2n\mu(X).
\]
Consequently,
\[
 \begin{aligned}
 0
 \leq\int_X|\nabla(\varphi-f)|^2\,d\mu
 =\int_X|\nabla f|^2\,d\mu-2n\mu(X)
 \leq0.
 \end{aligned}
\]
Thus \(\varphi-f\) is constant on
\(X_{\reg}\).
\end{proof}

\begin{proof}[Proof of Theorem \ref{cone}]
Since \(\varphi-f\) is constant on \(X_{\reg}\), we have
\[
 \sqrt{-1}\partial\bar\partial f=\omega,\qquad
 dV_g=e^{f-\varphi}\Omega=C\Omega
\]
for some \(C>0\).
Since \(\Omega\) has zero Ricci form,
\(\operatorname{Ric}_g=0\) on \(X_{\reg}\).
The shrinker equation and \(\nabla f(\varphi)=2f\) give
\[
 \nabla_g^2 f=g,\qquad
 |\nabla f|^2=2f.
\]
Along the flow \(\Phi_t\), we have
\begin{equation}\label{eq:homogeneity}
\frac{d}{dt}f(\Phi_t x)
=|\nabla f|^2(\Phi_t x)
=2f(\Phi_t x),\qquad f\circ\Phi_t=e^{2t}f,\qquad
 \Phi_t^*\omega=e^{2t}\omega.
\end{equation}
Set \(r=\sqrt{2f}\), then \(|\nabla r|=1\) and
\(\nabla f=r\nabla r\) on \(X_{\reg}\setminus\{o\}\).
Every orbit meets
\(\Sigma=\{r=1\}\cap X_{\reg}\) exactly once, and the map
\[
 (r,p)\longmapsto\Phi_{\log r}(p)
\]
identifies
\[
 (X_{\reg}\setminus\{o\},g)
 \cong
 \bigl((0,\infty)\times\Sigma,\,
       dr^2+r^2g_\Sigma\bigr).
\]

By continuity of $f$ and uniqueness of plurisubharmonic
extension, $f=\varphi+c$ holds on $X$.
Thus $f$ is nonnegative, locally bounded and plurisubharmonic
on $X$. 
Both $(\sqrt{-1}\partial\bar\partial f)^n$ and $\Omega$
give zero mass to $X\setminus X_{\reg}$.
Hence the equation on $X_{\reg}$ holds on all of $X$:
\[
 \frac{(\sqrt{-1}\partial\bar\partial f)^n}{n!}=C\Omega.
\]
Together with the cone structure on $X_{\reg}\setminus\{o\}$,
this shows that $\omega$ is a Ricci-flat K\"ahler cone metric
in the pluripotential sense; see
\cite[Section~2.5]{li2021notes}.
\end{proof}

\bibliography{ref}

@article{sun2024k,
  title={Kähler-Ricci shrinkers and Fano fibrations},
  author={Sun, Song and Zhang, Junsheng},
  journal={arXiv preprint arXiv:2410.09661},
  year={2024}
}

@article{conlon2024classification,
  title={Classification results for expanding and shrinking gradient Kähler--Ricci solitons},
  author={Conlon, Ronan J and Deruelle, Alix and Sun, Song},
  journal={Geometry \& Topology},
  volume={28},
  number={1},
  pages={267--351},
  year={2024},
  publisher={Mathematical Sciences Publishers}
}

@article{cifarelli2024finite,
  title={On finite time Type {I} singularities of the Kähler--Ricci flow on compact Kähler surfaces},
  author={Cifarelli, Charles and Conlon, Ronan J and Deruelle, Alix},
  journal={Journal of the European Mathematical Society},
  year={2024},
  publisher={European Mathematical Society}
}

@article{bamler2024new,
  title={A new complete two-dimensional shrinking gradient Kähler-Ricci soliton},
  author={Bamler, Richard H and Cifarelli, Charles and Conlon, Ronan J and Deruelle, Alix},
  journal={Geometric and Functional Analysis},
  volume={34},
  number={2},
  pages={377--392},
  year={2024},
  publisher={Springer}
}

@book{chow2025ricci,
  title={Ricci Solitons in Dimensions $4 $ and Higher},
  author={Chow, Bennett and Kotschwar, Brett and Munteanu, Ovidiu},
  volume={293},
  year={2025},
  publisher={American Mathematical Society}
}

@article{zhang2019gradient,
  title={Gradient Kähler--Ricci Solitons with Nonnegative Orthogonal Bisectional Curvature},
  author={Zhang, Shijin},
  journal={Results in Mathematics},
  volume={74},
  number={4},
  pages={127},
  year={2019},
  publisher={Springer}
}

@article{ni2005ancient,
  title={Ancient solutions to Kähler-Ricci flow},
  author={Ni, Lei},
  journal={Mathematical Research Letters},
  volume={12},
  number={5},
  pages={633--654},
  year={2005},
  publisher={International Press of Boston, Inc. Somerville, MA 02143, USA}
}

@article{munteanu2017positively,
  title={Positively curved shrinking Ricci solitons are compact},
  author={Munteanu, Ovidiu and Wang, Jiaping},
  journal={Journal of Differential Geometry},
  volume={106},
  number={3},
  pages={499--505},
  year={2017},
  publisher={Lehigh University}
}

@article{li2026Kähler,
  title={On K{\"a}hler Ricci shrinker surfaces},
  author={Li, Yu and Wang, Bing},
  journal={Acta Mathematica},
  volume={236},
  number={1},
  pages={1--50},
  year={2026},
  publisher={International Press of Boston, Inc. Somerville, MA 02143, USA}
}

@book{chow2023ricci,
  title={Ricci solitons in low dimensions},
  author={Chow, Bennett},
  volume={235},
  year={2023},
  publisher={American Mathematical Society}
}

@article{collins2019sasaki,
  title={Sasaki--Einstein metrics and K--stability},
  author={Collins, Tristan and Sz{\'e}kelyhidi, G{\'a}bor},
  journal={Geometry \& Topology},
  volume={23},
  number={3},
  pages={1339--1413},
  year={2019},
  publisher={Mathematical Sciences Publishers}
}

@article{chen2009strong,
  title={Strong uniqueness of the Ricci flow},
  author={Chen, Bing-Long},
  journal={Journal of Differential Geometry},
  volume={82},
  number={2},
  pages={363--382},
  year={2009},
  publisher={Lehigh University}
}

@incollection {MR954419,
    AUTHOR = {Hamilton, Richard S.},
     TITLE = {The {R}icci flow on surfaces},
 BOOKTITLE = {Mathematics and general relativity ({S}anta {C}ruz, {CA},
              1986)},
    SERIES = {Contemp. Math.},
    VOLUME = {71},
     PAGES = {237--262},
 PUBLISHER = {Amer. Math. Soc., Providence, RI},
      YEAR = {1988},
      ISBN = {0-8218-5079-2},
   MRCLASS = {53C20 (35K55 58G11)},
  MRNUMBER = {954419},
MRREVIEWER = {Dennis\ M.\ DeTurck},
       DOI = {10.1090/conm/071/954419},
       URL = {https://doi.org/10.1090/conm/071/954419},
}

@incollection {MR1375255,
    AUTHOR = {Hamilton, Richard S.},
     TITLE = {The formation of singularities in the {R}icci flow},
 BOOKTITLE = {Surveys in differential geometry, {V}ol.\ {II} ({C}ambridge,
              {MA}, 1993)},
     PAGES = {7--136},
 PUBLISHER = {Int. Press, Cambridge, MA},
      YEAR = {1995},
      ISBN = {1-57146-027-6},
   MRCLASS = {53C21 (58G30)},
  MRNUMBER = {1375255},
MRREVIEWER = {Man\ Chun\ Leung},
}

@article{naber2007noncompact,
  title={Noncompact shrinking four solitons with nonnegative curvature},
  author={Naber, Aaron},
  journal={Journal f{\"u}r die reine und angewandte Mathematik (Crelles Journal)},
  volume={2010},
  number={645},
  pages={125},
  year={2010},
  publisher={Walter de Gruyter GmbH}
}

@article{hamilton1982three,
  title={Three-manifolds with positive Ricci curvature},
  author={Hamilton, Richard S},
  journal={Journal of Differential geometry},
  volume={17},
  number={2},
  pages={255--306},
  year={1982},
  publisher={Lehigh University}
}

@article{LiZhang26,
  title={Gromov-Hausdorff Limits of Noncollapsed Kähler-Ricci Flows and the Geometry of Ricci Shrinkers},
  author={Li, Yu and Zhang, Junsheng},
  journal={arXiv preprint arXiv:2607.25644},
  year={2026}
}

@article{wu2016remarks,
  title={Remarks on shrinking gradient K{\"a}hler--Ricci solitons with positive bisectional curvature},
  author={Wu, Guoqiang and Zhang, Shijin},
  journal={Comptes Rendus. Math{\'e}matique},
  volume={354},
  number={7},
  pages={713--716},
  year={2016}
}

@misc{Perelman2003RicciFlowSurgery,
  author        = {Perelman, Grisha},
  title         = {Ricci flow with surgery on three-manifolds},
  year          = {2003},
  eprint        = {math/0303109},
  archivePrefix = {arXiv},
  primaryClass  = {math.DG},
  url           = {https://arxiv.org/abs/math/0303109}
}

@article{NiWallach2008Classification,
  author  = {Ni, Lei and Wallach, Nolan},
  title   = {On a classification of gradient shrinking solitons},
  journal = {Mathematical Research Letters},
  volume  = {15},
  number  = {5},
  pages   = {941--955},
  year    = {2008},
  doi     = {10.4310/MRL.2008.v15.n5.a9},
  eprint  = {0710.3194},
  archivePrefix = {arXiv},
  primaryClass  = {math.DG}
}

@article{wu2023shrinkers,
  title={Shrinkers with Curvature-Pinching Conditions are Compact: G. Wu, J.-Y. Wu},
  author={Wu, Guoqiang and Wu, Jia-Yong},
  journal={The Journal of Geometric Analysis},
  volume={33},
  number={9},
  pages={285},
  year={2023},
  publisher={Springer}
}

@article{LiNi2020Orthogonal,
  author  = {Li, Xiaolong and Ni, Lei},
  title   = {K{\"a}hler--Ricci shrinkers and ancient solutions with
             nonnegative orthogonal bisectional curvature},
  journal = {Journal de Math{\'e}matiques Pures et Appliqu{\'e}es},
  volume  = {138},
  pages   = {28--45},
  year    = {2020},
  doi     = {10.1016/j.matpur.2019.09.007}
}

@article{PetersenWylie2010Classification,
  author  = {Petersen, Peter and Wylie, William},
  title   = {On the classification of gradient Ricci solitons},
  journal = {Geometry \& Topology},
  volume  = {14},
  pages   = {2277--2300},
  year    = {2010},
  doi     = {10.2140/gt.2010.14.2277}
}

@article{LiNiWang2018PIC,
  author  = {Li, Xiaolong and Ni, Lei and Wang, Kui},
  title   = {Four-Dimensional Gradient Shrinking Solitons with
             Positive Isotropic Curvature},
  journal = {International Mathematics Research Notices},
  volume  = {2018},
  number  = {3},
  pages   = {949--959},
  year    = {2018},
  doi     = {10.1093/imrn/rnw269}
}

@article{QuWu2024Compact,
  author  = {Qu, Yuanyuan and Wu, Guoqiang},
  title   = {When are shrinking gradient Ricci soliton compact},
  journal = {Differential Geometry and its Applications},
  volume  = {93},
  pages   = {102102},
  year    = {2024},
  doi     = {10.1016/j.difgeo.2023.102102}
}

@article{berndtsson2009curvature,
  title={Curvature of vector bundles associated to holomorphic fibrations},
  author={Berndtsson, Bo},
  journal={Annals of mathematics},
  pages={531--560},
  year={2009},
  publisher={JSTOR}
}

@article{zhang2009completeness,
  title={On the completeness of gradient Ricci solitons},
  author={Zhang, Zhu-Hong},
  journal={Proceedings of the American Mathematical Society},
  volume={137},
  number={8},
  pages={2755--2759},
  year={2009}
}

@article{xu2026compactness,
  title={Compactness and Rigidity of Complete K$\backslash$" ahler Ricci Shrinkers},
  author={Xu, Tongxin and Zhang, Zhenlei},
  journal={arXiv preprint arXiv:2608.10953},
  year={2026}
}

@book{petersen2006riemannian,
  title={Riemannian geometry},
  author={Petersen, Peter},
  year={2006},
  publisher={Springer}
}

@article{HeinznerIannuzzi1997,
  title={Integration of local actions on holomorphic fiber spaces},
  author={Heinzner, Peter and Iannuzzi, Andrea},
  journal={Nagoya Mathematical Journal},
  volume={146},
  pages={31--53},
  year={1997},
  publisher={Cambridge University Press}
}

@article{PW09,
  author  = {Petersen, Peter and Wylie, William},
  title   = {Rigidity of Gradient Ricci Solitons},
  journal = {Pacific Journal of Mathematics},
  volume  = {241},
  number  = {2},
  pages   = {329--345},
  year    = {2009},
  doi     = {10.2140/pjm.2009.241.329}
}

@article{FLGR16,
  author  = {Fern{\'a}ndez-L{\'o}pez, Manuel
             and Garc{\'i}a-R{\'i}o, Eduardo},
  title   = {On Gradient Ricci Solitons with Constant
             Scalar Curvature},
  journal = {Proceedings of the American Mathematical Society},
  volume  = {144},
  number  = {1},
  pages   = {369--378},
  year    = {2016},
  doi     = {10.1090/proc/12693}
}

@article{CZ23,
  author  = {Cheng, Xu and Zhou, Detang},
  title   = {Rigidity of Four-Dimensional Gradient
             Shrinking Ricci Solitons},
  journal = {Journal f{\"u}r die reine und angewandte Mathematik},
  volume  = {2023},
  number  = {802},
  pages   = {255--274},
  year    = {2023},
  doi     = {10.1515/crelle-2023-0042}
}

@article{LOQW24,
  author        = {Li, Fengjiang and Ou, Jianyu
                   and Qu, Yuanyuan and Wu, Guoqiang},
  title         = {Rigidity of Five-Dimensional Shrinking
                   Gradient Ricci Solitons with Constant
                   Scalar Curvature},
  year          = {2024},
  journal        = {arxiv:2411.10712},
}

@article{LOQW25,
  author        = {Li, Fengjiang and Ou, Jianyu
                   and Qu, Yuanyuan and Wu, Guoqiang},
  title         = {Rigidity of Five-Dimensional Shrinking
                   Gradient Ricci Solitons},
  year          = {2025},
  journal       = {arxiv:2506.00887},
}

@article{HLOWZ26,
  author        = {Hor{\'a}cio, Matheus and Li, Wenqi
                   and Ou, Jianyu and Wu, Guoqiang
                   and Zhou, Detang},
  title         = {Rigidity of {K\"ahler}--Ricci Solitons
                   with Constant Scalar Curvature},
  year          = {2026},
  journal={arXiv:2610.11935}
}

@article{WW26,
  author        = {Wang, Chen and Wu, Guoqiang},
  title         = {A Note on Rigidity of Shrinking Gradient
                   Ricci Solitons with Constant Scalar Curvature},
  year          = {2026},
  journal        = {arXiv:2604.23939}
}

@article{LQW26,
  author        = {Li, Fengjiang and Qu, Yuanyuan
                   and Wu, Guoqiang},
  title         = {Rigidity of Shrinking Gradient Ricci Soliton
                   with Constant Scalar Curvature},
  year          = {2026},
  journal        = {arxiv:2608.20040},
}

@article{SunZhangII,
  author = {Sun, Song and Zhang, Junsheng},
  title  = {{K\"ahler--Ricci} shrinkers and {Fano} fibrations, {II}},
  note   = {Preprint},
  journal    = {https://zhang-junsheng.github.io/assets/kahler-ricci-shrinkers-and-fano-fibrations-II.pdf}
}

@article{SJ,
  title={Riemannian geometry of Kahler-Einstein currents},
  author={Song, Jian},
  journal={arXiv preprint arXiv:1404.0445},
  year={2014}
}

@article{li2021notes,
  title={Notes on weighted K{\"a}hler-Ricci solitons and application to Ricci-flat K{\"a}hler cone metrics},
  author={Li, Chi},
  journal={arXiv preprint arXiv:2107.02088},
  volume={1},
  year={2021}
}

@article{ConlonDeruelle26,
  title={Uniqueness of shrinking Kähler-Ricci solitons on resolutions of Kähler cones},
  author={Conlon, Ronan J and Deruelle, Alix},
  journal={arXiv preprint arXiv:2607.26054},
  year={2026}
}
\bibliographystyle{alpha}

\bigskip

\noindent
\textsc{Tongxin Xu}\\
School of Mathematical Sciences, Capital Normal University, Beijing, 100048, China\\
\textit{Email address:} \texttt{2250501032@cnu.edu.cn}

\bigskip

\noindent
\textsc{Zhenlei Zhang}\\
School of Mathematical Sciences, Capital Normal University, Beijing, 100048, China\\
\textit{Email address:} \texttt{zhleigo@aliyun.com}
\end{document}